\documentclass{amsart}
\usepackage{amssymb}

\usepackage{amsfonts}

\begin{document}
\title[A remark on the trigonometric system]{A remark on the trigonometric
system}
\author{Alexander Kushpel}
\address{Department of Mathematics\\
University of Leicester}
\email{ak412@le.ac.uk}
\date{May 25, 2015}

\begin{abstract}
We present a new property of the trigonometric system arranged in a natural
order. It is shown that the sequence of subspaces of trigonometric
polynomials $\mathcal{T}_{n}$ is optimal in the sense of order of
convergence on convolution classes $K\ast U_{p}$ in $L_{q}$ for any $%
1<p,q<\infty $ just in the cases of "very slow" or "very fast" rate of decay
of Fourier coefficients of $K$.
\end{abstract}

\maketitle

\section{Introduction}

Let $X$ be a Banach space and $A$ a convex, compact, centrally symmetric
subset of $X$. The Kolmogorov $n$-width of $A$ in $X$ is defined by%
\[
d_{n}\left( A,X\right) :=\inf_{X_{n}\subset X}\sup_{f\in A}\inf_{g\in
X_{n}}\left\Vert f-g\right\Vert _{X},
\]
where $X_{n}$ runs over all subspaces of $X$ of dimension $\leq n$. Let $%
\mathbb{T}^{1}$ be the unit circle with the Lebesgue measure $dx$. The space
of $p$-integrable functions $\phi $, i.e. such that $\left\Vert \phi
\right\Vert _{p}<\infty $ is denoted by $L_{p}\left( \mathbb{T}^{1}\right)
\equiv L_{p}$, where
\[
\left\Vert \phi \right\Vert _{p}:=\left( \int_{\mathbb{T}^{1}}\left\vert
\phi \left( x\right) \right\vert ^{p}dx\right) ^{1/p},\text{ \ }1 \leq
p<\infty .
\]
Consider the sequence of subspaces $\mathcal{T}_{n}$ of trigonometric
polynomials with the usual order,
\[
\mathcal{T}_{0}:=\mathrm{lin}\left\{ 1\right\} ,\cdot \cdot \cdot ,\mathcal{T%
}_{n}:=\mathrm{lin}\left\{ 1,\cos kx,\sin kx,1\leq k\leq n\right\} ,\cdot
\cdot \cdot .
\]

Let $\phi \in L_{p}$ with the formal Fourier series
\[
\phi \sim \sum_{k=1}^{\infty }a_{k}\left( \phi \right) \cos kx+b_{k}\left(
\phi \right) \sin kx.
\]%
A wide range of sets of smooth functions on $\mathbb{T}^{1}$ can be
introduced using multipliers $\Lambda :=\left\{ \lambda _{k},\text{ }k\in
\mathbb{N}\right\} $ (see, e.g. \cite{step}). We say that $f\in \Lambda
_{\beta }U_{p},$ $\beta \in \mathbb{R}$ if
\[
f\sim \sum_{k=1}^{\infty }\lambda _{k}\left( a_{k}\left( \phi \right) \cos
\left( kx-\frac{\beta \pi }{2}\right) +b_{k}\left( \phi \right) \sin \left(
kx-\frac{\beta \pi }{2}\right) \right) ,
\]%
where $\phi \in U_{p}:=\left\{ \phi \left\vert \left\Vert \phi \right\Vert
_{p}\leq 1\right. \right\} .$ If there exists a function $K\in L_{1},$%
\[
K\sim \sum_{k=1}^{\infty }\lambda _{k}\cos \left( kx-\frac{\beta \pi }{2}%
\right) ,
\]%
then the set $\Lambda _{\beta }U_{p}$ can be represented in the convolution
form
\[
f\left( x\right) =\frac{1}{2\pi }\int_{\mathbb{T}^{1}}K\left( x-y\right)
\phi \left( y\right) dy.
\]%
In this case we say that $f\in K\ast U_{p}$. If $K\notin L_{1}$ then we
should consider generalised convolutions. In particular, if $\lambda
_{k}=k^{-r},$ $\beta =r,$ $r>0$ we get Sobolev's classes $W_{p}^{r}$. If $%
\lambda _{k}=\exp \left( -\mu k^{r}\right) ,$ $\beta =0,$ $\mu >0,$ \ $0<r<1$%
, then the class $K\ast U_{p}$ consists of infinitely differentiable
functions. In the case $\lambda _{k}=\exp \left( -\mu k^{r}\right) ,$ $\beta
=0,$ $\mu >0,$ \ $r=1$ $\left( r>1\right) $ we get classes of analytic
(entire) functions respectively.

Let $\lambda _{k}=k^{-r},$ $\beta =r.$ It is well-known that in this case
("finite smoothness", i.e. $r>1/p$ if $1<p\leq 2\leq q<\infty $, or $2\leq
p\leq q<\infty $)
\[
d_{n}\left( W_{p}^{r},L_{q}\right) \asymp \left\{
\begin{array}{cc}
n^{-r},r>0, & 1<q<p<\infty , \\
n^{-r+1/p-1/q},r>1/p-1/q, & 1<p\leq q\leq 2, \\
n^{-r+1/p-1/2},r>1/p, & 1<p\leq 2\leq q<\infty , \\
n^{-r},r>1/p, & 2\leq p\leq q<\infty ,%
\end{array}%
\right.
\]
as $n\rightarrow \infty $. It is easy to show that
\[
E_{n}\left( W_{p}^{r},L_{q}\right) :=\sup_{f\in W_{p}^{r}}\inf_{g\in
\mathcal{T}_{n}}\left\Vert f-g\right\Vert _{L_{q}}\asymp n^{-r+\left(
1/p+1/q\right) _{+}},\text{ \ }r>1/p-1/q,
\]
where $\left( a\right) _{+}:=\max \left\{ a,0\right\} $. Hence the sequence
of trigonometric polynomials $\mathcal{T}_{n}$ is optimal if $r>1/p-1/q,$ $%
1<p\leq q\leq 2$, or $r>0,$ $1<q<p<\infty $. In the remaining cases, $r>1/p,$
$1<p\leq 2\leq q<\infty $, $r>1/p,$ $2\leq p\leq q<\infty ,$ $\mathcal{T}%
_{n} $ is not optimal. If $1/p-1/q<r<1/p$ then we have the case of so-called
"small smoothness",%
\[
d_{n}\left( W_{p}^{r},L_{q}\right) \asymp n^{\gamma },1/p-1/q<r<1/p,r\neq
\beta ,2\leq p\leq q<\infty ,
\]
where
\[
\gamma =\max \left\{ -r,q\left( -r+1/p-1/q\right) /2\right\}
\]
and%
\[
\beta =\frac{1/p-1/q}{2\left( 1/2-1/q\right) }.
\]
In the case
\[
1/p-1/q<r<1/p,1<p<2<q<\infty
\]
we have%
\[
d_{n}\left( W_{p}^{r},L_{q}\right) \asymp n^{q\left( -r+1/p-1/q\right) /2}.
\]
This means that the sequence of subspaces $\mathcal{T}_{n}$ is not optimal
in this case. If $\lambda _{k}=\exp \left( -\mu k^{r}\right) ,$ $\beta =0,$ $%
\mu >0,$ \ $0<r<1$, then
\[
d_{2n}\left( W_{p}^{r},L_{q}\right) \asymp \left\{
\begin{array}{cc}
\exp \left( -\mu n^{r}\right) n^{\left( 1-r\right) \left( 1/p-1/q\right) },
& 1<p\leq q\leq 2, \\
\exp \left( -\mu n^{r}\right) , & 1<q\leq p\leq 2\text{ or }2\leq p,q<\infty
, \\
\exp \left( -\mu n^{r}\right) n^{\left( 1-r\right) \left( 1/p-1/2\right) } &
1<p\leq 2\leq q<\infty .%
\end{array}%
\right.
\]
It is easy to show that
\[
E_{2n}\left( K\ast U_{p},L_{q}\right) :=\sup_{f\in K\ast U_{p}}\inf_{g\in
\mathcal{T}_{n}}\left\Vert f-g\right\Vert _{L_{q}}
\]
\[
\asymp \exp \left( -\mu n^{r}\right) n^{\left( 1-r\right) _{+}\left(
1/p-1/q\right) },
\]
$1<p,q<\infty $, where
\[
K(x)=\sum_{k=1}^{\infty }\exp \left( -\mu k^{r}\right) \cos kx.
\]
Hence $\mathcal{T}_{n}$ is not optimal if $2\leq p,q<\infty $ and $1<p\leq
2\leq q<\infty $.

Finally, in the case of analytic or entire smoothness, i.e. if $r\geq 1$ we
have%
\[
d_{2n}\left( W_{p}^{r},L_{q}\right) \asymp E_{2n}\left( K\ast
U_{p},L_{q}\right) \asymp \exp \left( -\mu n^{r}\right) .
\]
This means the the sequence of subspaces $\mathcal{T}_{n}$ is optimal for $%
K\ast U_{p}$ in $L_{q}$ for any $1<p,q<\infty $. See \cite{kushpel1} for
more information.

\section{The results}

In this section we prove the following statement.

\textbf{Theorem 1.} 
\emph{Let $\lambda _{k}=k^{-\left( 1/p-1/q\right) _{+}}\left( \ln k\right)
^{-\gamma },$ $\gamma >0,$ $1<p,q<\infty $, then the sequence of subspaces $%
\mathcal{T}_{n}$ is optimal in the sense of order as $n\rightarrow \infty $.
}

\textbf{Proof} Let $\lambda _{k}=k^{-\left( 1/p-1/q\right) _{+}}\left( \ln
k\right) ^{-\gamma },$ $\gamma >0,$ $1<p,q<\infty $. Using results from \cite%
{kushpel0, kushpel2} it is possible to show that
\[
d_{n}\left( \Lambda U_{p},L_{q}\right) \leq E_{n}\left( K\ast
U_{p},L_{q}\right) \ll \left( \ln k\right) ^{-\gamma },\text{ }1<p<q<\infty
.
\]
We turn to the lower bounds now. Clearly,
\[
d_{n}\left( \Lambda U_{p},L_{q}\right) \geq d_{n}\left( \Lambda U_{p}\cap
\mathcal{T}_{m},L_{q}\right)
\]
and%
\[
d_{n}\left( \Lambda U_{p}\cap \mathcal{T}_{m},L_{q}\right) \geq
C_{q}^{-1}d_{n}\left( \Lambda U_{p}\cap \mathcal{T}_{m},L_{q}\cap \mathcal{T}%
_{m}\right)
\]
since $\left\Vert S_{m}:L_{q}\rightarrow L_{q}\cap \mathcal{T}%
_{m}\right\Vert \leq C_{q},$ $1<q<\infty $, where $S_{m}$ is the operator of
orthogonal projection onto $\mathcal{T}_{m}$. Finally, applying
Marcinkiewicz-Zygmund inequality%
\[
C_{1,p}\left\Vert t_{m}\right\Vert _{p}\leq \left(
m^{-1}\sum_{k=1}^{2m+1}\left\vert t_{n}\left( \frac{2\pi k}{2m+1}\right)
\right\vert ^{p}\right) ^{1/p}\leq C_{2,p}\left\Vert t_{m}\right\Vert _{p},%
\text{ \ }\forall t_{m}\in \mathcal{T}_{m}
\]
we get
\[
d_{n}\left( \Lambda U_{p},L_{q}\right) \gg \left( \ln m\right) ^{-\gamma
}d_{n}\left( B_{p}^{m},l_{q}^{m}\right) ,
\]
where the norm in $l_{q}^{m},1\leq q\leq \infty $ is defined as usual, $%
\left\Vert v\right\Vert :=\left( \sum_{k=1}^{m}\left\vert v_{k}\right\vert
^{q}\right) ^{1/q}$ and $B_{p}^{m}$ is the unit ball in $l_{p}^{m}$. We will
need the following result \cite{kashin1, kashin2}%
\[
C_{1}\left( p,q\right) \leq \frac{d_{n}\left( B_{p}^{m},l_{q}^{m}\right) }{%
\Phi \left( m,n,p,q\right) }\leq C_{2}\left( p,q\right)
\]
for any $m>n$, where%
\[
\Phi \left( m,n,p,q\right) :=\left( \min \left\{ 1,m^{1/q}n^{-1/2}\right\}
\right) ^{\left( 1/p-1/q\right) /\left( 1/2-1/q\right) },
\]
if $2\leq p<q\leq \infty $ and%
\[
\Phi \left( m,n,p,q\right) :=\max \left\{ m^{1/q-1/p},\min \left\{
1,m^{1/q}n^{-1/2}\right\} \left( 1-n/m\right) ^{1/2}\right\}
\]
if $1\leq p<2\leq q\leq \infty $. Let $2\leq p<q\leq \infty $ then
\[
d_{n}\left( \Lambda U_{p},L_{q}\right) \gg \left( \ln m\right) ^{-\gamma
}d_{n}\left( B_{p}^{m},l_{q}^{m}\right)
\]
\[
\gg \left( \ln m\right) ^{-\gamma }\left( \min \left\{
1,m^{1/q}n^{-1/2}\right\} \right) ^{\left( 1/p-1/q\right) /\left(
1/2-1/q\right) }.
\]
Let $m=n^{q/2}$ then
\[
d_{n}\left( \Lambda U_{p},L_{q}\right) \gg \left( \ln m\right) ^{-\gamma
}\gg \left( \ln n\right) ^{-\gamma }.
\]
If $1<p<2\leq q<\infty $ then%
\[
d_{n}\left( \Lambda U_{p},L_{q}\right) \gg \left( \ln m\right) ^{-\gamma
}\left( \min \left\{ 1,m^{1/q}n^{-1/2}\right\} \right) ^{\left(
1/p-1/q\right) /\left( 1/2-1/q\right) }
\]
\[
\gg \left( \ln m\right) ^{-\gamma }.
\]
Hence the trigonometric system is an optimal in the logarithmic
neighbourhood of embedding of $L_{p}$ into $L_{q}$.\,\,\,  

\end{document}